\documentclass[11pt,a4paper,final]{article}%
\usepackage{amscd}
\usepackage{amsmath,amsthm}
\usepackage{graphicx}
\usepackage{amsfonts}
\usepackage{amssymb}
\usepackage{showkeys}%
\newtheorem{theorem}{Theorem}[section]
\newtheorem{lemma}[theorem]{Lemma}
\newtheorem{corollary}[theorem]{Corollary}

\newtheorem{definition}[theorem]{Definition}
\newtheorem{remark}[theorem]{Remark}

\numberwithin{equation}{section}

\begin{document}

\title{On the existence of infinite energy solutions for
nonlinear Schr\"{o}dinger equations}
\author{Pablo Braz e Silva\footnotemark[1] \\{\small Departamento de Matem\'atica }\\
{\small Universidade Federal de Pernambuco }
\\{\small CEP 50740-540, Recife-PE, Brazil.}\\{\small \texttt{email:\ pablo@dmat.ufpe.br}}
\\ \\ Lucas C. F. Ferreira \\ {\small Departamento de Matem\'atica}
\\{\small Universidade Federal de Pernambuco}
\\{\small CEP 50740-540, Recife-PE, Brazil.}\\{\small \texttt{email:\ lcff@dmat.ufpe.br}}
\\ \\E. J. Villamizar-Roa
\\{\small Universidad Industrial de Santander}
\\{\small {\ A.A. 678. Bucaramanga, Colombia.}}
\\{\small \texttt{email:\ jvillami@uis.edu.co}}}
\date{}
\maketitle

\renewcommand{\thefootnote}{\fnsymbol{footnote}}

\footnotetext[1]{P. Braz e Silva is partly supported by
CAPES/MECD-DGU Brazil/Spain, grant \#117/06.}

\thanks{Math Subject Classification: 35Q55, 35D05, 35B40}


\renewcommand{\thefootnote}{\arabic{footnote}}

\begin{abstract}
We derive new results about existence and uniqueness of local and
global solutions for nonlinear Schr\"{o}dinger equation, including
self-similar global solutions. Our analysis is performed in the
framework of Marcinkiewicz spaces.

\end{abstract}

\section{ Introduction}

We consider the nonlinear Schr\"{o}dinger equation
\begin{align}
i\partial_{t}u+\Delta u  &  =\lambda|u|^{\rho}u,\ x\in\mathbb{R}^{n}%
,\ t\in\mathbb{R},\label{sc1}\\
u(0,x)  &  =\phi(x),\ x\in\mathbb{R}^{n}, \label{sc2}%
\end{align}
where $u=u(t,x)$ is a complex valued function, $\lambda$ is a fixed
complex number, and $0<\rho<\infty$. The initial value $\phi:\mathbb{R}^{n}%
\rightarrow\mathbb{C}$ is given. The Cauchy problem (\ref{sc1})-(\ref{sc2}) is
formally equivalent to the integral equation
\begin{equation}
u(t)=S(t)\phi-i\lambda\int_{0}^{t}S(t-s)(|u(s)|^{\rho}u(s))ds, \label{sc6}%
\end{equation}
where $S(t)$ is the unitary group determined by the linear Schr\"{o}dinger
equation
\[
\partial_{t}u-i\Delta u=0,\ x\in\mathbb{R}^{n},\ t\in\mathbb{R}.
\]
If $\phi\in\mathcal{S}(\mathbb{R}^{n})$ and $u$ is defined by $\widehat
{u(t)}(\xi)=e^{-i|\xi|^{2}t}\widehat{\phi}(\xi),$ for $\xi\in\mathbb{R}^{n}$,
then
\[
\widehat{u}_{t}+i|\xi|^{2}\widehat{u}=0
\]
in $\mathbb{R}\times\mathbb{R}^{n}$. In this case, the solution of
\begin{align}
\partial_{t}u-i\Delta u  &  =0,\ x\in\mathbb{R}^{n},\ t\in\mathbb{R}%
,\label{sc4}\\
u(0,x)  &  =\phi(x),\ x\in\mathbb{R}^{n}, \label{linearsh}%
\end{align}
is given by $u(t)=S(t)\phi=K_{t}\ast\phi$, where $K_{t}(x)=\left(
e^{-i|\xi|^{2}t}\right)  ^{\check{}}.$

Existence and uniqueness of local and global solutions of problem
(\ref{sc1})-(\ref{sc2}) have been much studied in the framework of the Sobolev
spaces $H^{s}$, $s\geq0$, i.e, the solutions and their derivatives have finite
energy. See, for instance, Ginibre and Velo \cite{Ginibre}-\cite{Ginibre3},
Kato \cite{Ka1}-\cite{Ka2}, Cazenave and Weissler\cite{Cazenave1}
-\cite{Cazenave4}, and the references therein.

As far as we know, the first authors to study infinite energy solutions of
(\ref{sc1})-(\ref{sc2}) were Cazenave and Weissler in \cite{Cazenave5}. There,
they consider the space
\[
\displaystyle X_{\rho}=\{u\in L_{loc}^{\infty}((0,\infty),\,L^{\rho
+2}(\mathbb{R}^{n}));\;\sup_{t>0}t^{\frac{\alpha}{2}}\Vert u(t)\Vert
_{L^{\rho+2}}<\infty\},
\]
where $\frac{\alpha}{2}=\frac{1}{\rho}-\frac{n}{2(\rho+2)}$ and
$\displaystyle\Vert\cdot\Vert_{L^{\rho+2}}$ denotes the usual $L^{\rho+2}$
norm. Under a suitable smallness condition on the initial data, they prove the
existence of global solutions of (\ref{sc1})-(\ref{sc2}) in $X_{\rho}$, for
$\rho$ in the range
\begin{equation}
\frac{\rho+2}{\rho+1}<\frac{n\rho}{2}<\rho+2. \label{range1}%
\end{equation}
If $n=1$ or $n=2$, condition (\ref{range1}) is equivalent to $\rho_{0}%
<\rho<\infty$, where $\rho_{0}$ is the positive value of $\rho$ for which
$\frac{\rho+2}{\rho+1}=\frac{n\rho}{2}$. If $n \geq3$, it is equivalent to
$\rho_{0}<\rho<\frac{4}{n-2}$. Later on, in \cite{Cazenave-Vega}, the Cauchy
problem (\ref{sc1})-(\ref{sc2}) was studied in the framework of weak-$L^{p}$
spaces. Using a Strichartz-type inequality, the authors obtained existence of
solutions in the class $L^{(p,\infty)}(\mathbb{R}^{n+1})$ $\equiv
L_{t}^{(p,\infty)}\left(  L_{x}^{(p,\infty)}\right)  ,$ where $(t,x)\in
\mathbb{R}\times\mathbb{R}^{n}$ and $p=\frac{\rho(n+2)}{2(\rho+1)}$, for
$\rho$ in the range
\begin{equation}
\rho_{0}<\frac{4(n+1)}{n(n+2)}<\rho<\frac{4(n+1)}{n^{2}}<\frac{4}{n-2}.
\label{range2}%
\end{equation}
In \cite{Planchon1}, the existence of solutions with initial data in
the Besov space $\dot{B}_{2}^{s_{\rho},\infty}(\mathbb{R}^{n})$,
with positive regularity $s_{\rho}=\frac{n}{2}-\frac{2}{\rho}>0$,
was proved for $\rho$ in the range
$\rho_{0}<\frac{4}{n}<\rho<\infty$. Note that if
$f\in\dot{B}_{2}^{s_{\rho },\infty}(\mathbb{R}^{n})$, then $f$ has
at least local finite energy.

We study equation (\ref{sc6}) in functional spaces of infinite energy. In the
first theorem proved here, we consider the initial data $\phi$ belonging to
the Marcinkiewicz space $L^{(\frac{\rho+2}{\rho+1},\infty)}$, and show
existence and uniqueness of local in time solutions in the class
\[
\displaystyle E_{\alpha,\beta}^{T}=\big\{u;\;\Vert u\Vert_{\alpha,\beta}%
=\sup_{-T<t<T}\left\vert t\right\vert ^{\frac{\alpha-\beta}{2}}\Vert
u(t)\Vert_{(\rho+2,\infty)}<\infty\big\},
\]
where $\frac{(\alpha-\beta)}{2}=\frac{n\rho}{2(\rho+2)}$, with $\frac{n\rho
}{2}<\frac{\rho+2}{\rho+1}$.
%Note that the norm $\Vert\cdot\Vert_{\alpha,\beta}$
%is not invariant by the scaling $u_{\mu}(t,x)=\mu^{\frac{2}{\mu}}u(\mu
%^{2}t,\mu x).$ This is a key point in order to obtain local in time solutions
%in the framework of Marcinkiewicz spaces.
Note that $\frac{n\rho}{2}<\frac{\rho+2}{\rho+1}$ is equivalent to
$0<\rho<\rho_{0}<\frac{4}{n}.$
%(note that$\rho_{0}<\frac{4}{n})$.
So, our range for $\rho$ is different from the ones in
\cite{Cazenave5, Cazenave-Vega, Planchon1}. The norm
$\Vert\cdot\Vert_{\alpha,\beta}$ is not invariant by the scaling
$u_{\mu}(t,x)=\mu^{\frac{2}{\rho}}u(\mu ^{2}t,\mu x).$ This is a key
point to obtain local in time solutions in Marcinkiewicz spaces. It
is also worth noting that our result allows one to consider singular
initial data as, for example, homogeneous functions $\left\vert
x\right\vert ^{-\frac{n(\rho+1)}{\rho+2}}\in L^{(\frac{\rho
+2}{\rho+1},\infty)}$.

Our second theorem concerns global in time solutions. We show
existence of such solutions in norms of type
$\sup_{|t|>0}|t|^{\alpha/2}\Vert u(t)\Vert _{L^{(\rho+2,\infty)}}$,
where $\frac{\alpha}{2}=\frac{1}{\rho}-\frac {n}{2(\rho+2)}$ and
\begin{equation}
\label{range3}\rho_{0}<\rho<\frac{4}{n-2}.
\end{equation}
This extends the result of Cazenave and Weissler \cite{Cazenave5} to the
context of Lorentz spaces.
%More
%specifically, we obtain existence of global solutions in norms of the type
%$\sup_{t>0}|t|^{\alpha/2}\Vert u(t)\Vert_{L^{(\rho+2,\infty)}},$ where
%$\frac{\alpha}{2}=\frac{1}{\rho}-\frac{n}{2(\rho+2)}$ and $\rho_{0}<\rho
%<\frac{4}{n-2}$.
Note that range (\ref{range3}) is greater than range (\ref{range2}).

As a corollary, we show that when the initial data $\phi$ is a homogeneous
function of degree $-\frac{2}{\rho}$, we obtain a self-similar solution, if
$\Vert S(1)\phi\Vert_{(\rho+2,\infty)}$ is sufficiently small. Moreover, we
discuss asymptotic stability of the global solutions, and show that regular
perturbations of the linear Schr\"{o}dinger equations are negligible for large
times. We also analyze the behavior of the local solutions as $t\rightarrow0$
in the space $L^{(\rho+2,\infty)}.$

Our approach is different from the methods used in
\cite{Cazenave-Vega,Planchon1}, where the authors use a Strichartz-type
inequality in weak-$L^{p}$ and Besov spaces, respectively. Indeed, our
existence results are based on bounds for the Schr\"{o}dinger linear group
$S(t)$ in the context of Lorentz spaces. In Lemma \ref{estim1}, we state and
prove these bounds via real interpolation techniques. They generalize the
bounds for usual $L^{p}$ spaces used in \cite{Cazenave5}.

In section \ref{results}, we carefully state our results and discuss their
improvement in the light of previous results. We prove them in section
\ref{proofs}.

\section{Main Results}

\label{results}

We first recall some facts about the Lorentz spaces. For more details see, for
instance, \cite{BL} and \cite{SW}.

Let $1<p\leq\infty$ and $1\leq q\leq\infty.$ A measurable function $f$ defined
on $\mathbb{R}^{n}$ belongs to Lorentz space $L^{(p,q)}(\mathbb{R}^{n})$ if
the quantity
\[
\Vert f\Vert_{(p,q)}=\left\{
\begin{array}
[c]{lll}%
\displaystyle \left(  \frac{p}{q}\int_{0}^{\infty}\left[  t^{\frac{1}{p}%
}f^{\ast\ast}(t)\right]  ^{q}\frac{dt}{t}\right)  ^{\frac{1}{q}} &
\mbox{, if } & 1<p<\infty\mbox{, }1\leq q<\infty, \vspace{.3cm}\\
\displaystyle \sup_{t>0}t^{\frac{1}{p}}f^{\ast\ast}(t) & \mbox{, if } &
1<p\leq\infty\mbox{, }q=\infty,
\end{array}
\right.
\]
is finite, where the $f^{\ast\ast}$ is defined for $t>0$ by
\begin{align*}
f^{\ast\ast}(t)  &  =  \frac{1}{t}\int_{0}^{t}f^{\ast}(s)\mbox{ } ds ,
\end{align*}
where
\begin{align*}
f^{\ast}(t)  &  =  \inf\big\{s>0 ; m\{{x}\in\mathbb{R}^{n}:|f({x})|>s\}\leq
t\big\}, \, t >0 .
\end{align*}
Note that $L^{p}(\mathbb{R}^{n})=L^{(p,p)}(\mathbb{R}^{n})$. The spaces
$L^{(p,\infty)}(\mathbb{R}^{n})$ are called weak-$L^{p}$ spaces or
Marcinkiewicz spaces. Lorentz spaces have the same scaling relation as $L^{p}$
spaces, that is, for all $\lambda>0$ one has $\Vert f(\lambda x)\Vert
_{(p,q)}=\lambda^{-\frac{n}{p}}\Vert f\Vert_{(p,q)},$ where $1\leq p<\infty$
and $1\leq q\leq\infty$. Morover, Lorentz spaces can be constructed via real
interpolation \cite{BL}. Indeed,
\[
L^{(p,q)}(\mathbb{R}^{n})=(L^{1}(\mathbb{R}^{n}),L^{\infty}(\mathbb{R}%
^{n}))_{1-\frac{1}{p},q},\ 1<p<\infty.
\]
They have the interpolation property
\[
(L^{(p_{0},q_{0})}(\mathbb{R}^{n}),L^{(p_{1},q_{1})}(\mathbb{R}^{n}%
))_{\theta,q}=L^{(p,q)}(\mathbb{R}^{n}),
\]
provided $0<p_{0}<p_{1}<\infty,\ 0<\theta<1,\
\frac{1}{p}=\frac{1-\theta }{p_{0}}+\frac{\theta}{p_{1}},\ 1\leq
q_{0},q_{1},q\leq\infty$, where $(\cdot,\cdot)_{\theta,q}$ stands
for the real interpolation spaces constructed via the $K$-method
\cite{BL}.

We begin by bounding the Schr\"{o}dinger group $S(t)$\ in Lorentz spaces.

\begin{lemma}
\label{estim1} Let $1\leq d\leq\infty$, and $1<p<2$. If $p^{\prime}$ is such
that $\frac{1}{p} + \frac{1}{p^{\prime}} = 1 $, then there exists a constant
$C=C(n,\gamma,p)>0$ such that
\begin{equation}
\left\Vert S(t)\varphi\right\Vert _{(p^{\prime},d)}\leq C|t|^{-\frac{n}%
{2}(\frac{2}{p}-1)}\left\Vert \varphi\right\Vert _{(p,d)}, \label{semigroup}%
\end{equation}
for all $\varphi\in L^{(p,d)}(\mathbb{R}^{n} )$ and all $t\neq0.$
\end{lemma}

\begin{proof}
Fix $t\neq0$ and let $1<p_{0}<p<p_{1}<2$ such that $\frac{1}{p^{\prime}}%
=\frac{\lambda}{p_{0}}+\frac{1-\lambda}{p_{1}}$ and $0<\lambda<1$. By the well
known $L^{p}=L^{(p,p)}$ estimate of Schr\"{o}dinger group, we have that
$
S(t):L^{p_{0}}\rightarrow L^{p_{0}^{\prime}}$ and $S(t):L^{p_{1}%
}\rightarrow L^{p_{1}^{\prime}}$,
where the operator norms are respectively bounded by
\begin{eqnarray*}
\left\Vert S(t)\right\Vert _{p_{0}\rightarrow p_{0}^{\prime}} & \leq &
C|t|^{-\frac{n}{2}(\frac{2}{p_{0}}-1)},
\\
\left\Vert S(t)\right\Vert
_{p_{1}\rightarrow p_{1}^{\prime}} & \leq & C|t|^{-\frac{n}{2}(\frac{2}{p_{1}}%
-1)}.
\end{eqnarray*}
Through real interpolation,
\begin{align*}
\left\Vert S(t)\right\Vert _{(p,d)\rightarrow(p^{\prime},d)}  &
\leq\left\Vert S(t)\right\Vert _{p_{0}\rightarrow p_{0}^{\prime}}^{\lambda
}\left\Vert S(t)\right\Vert _{p_{1}\rightarrow p_{1}^{\prime}}^{1-\lambda}\\
&  \leq\left(  C|t|^{-\frac{n}{2}(\frac{2}{p_{0}}-1)}\right)  ^{\lambda
}\left(  C|t|^{-\frac{n}{2}(\frac{2}{p_{1}}-1)}\right)  ^{1-\lambda}\\
&  =C|t|^{-\frac{n}{2}(\frac{2}{p}-1)},
\end{align*}
which is equivalent to (\ref{semigroup}).
\end{proof}
From now on, we denote $\displaystyle \alpha:=\frac{2}{\rho}-\frac{n}{\rho+2}$
and $\displaystyle \beta:=\frac{2}{\rho}-\frac{n(\rho+1)}{(\rho+2)}$.

\begin{definition}
Let $0<\rho<\infty$ and $0<T\leq\infty$. We denote by $E_{\alpha}$ and
$E_{\alpha,\beta}^{T}$ the Banach spaces
\begin{align}
E_{\alpha}  &  = \left\{  u ; \left\vert t\right\vert ^{\frac{\alpha}{2}}u\in
BC((-\infty,\infty);L^{(\rho+2,\infty)})\right\} ,\\
E_{\alpha,\beta}^{T}  &  =\left\{  u ; \left\vert t\right\vert ^{\frac
{\alpha-\beta}{2}}u\in BC((-T,T);L^{(\rho+2,\infty)})\right\} ,
\end{align}
with respective norms
\[
\Vert u\Vert_{\alpha}=\sup_{-\infty<t<\infty}\left\vert t\right\vert
^{\frac{\alpha}{2}}\Vert u(t)\Vert_{(\rho+2,\infty)},
\]
and
\[
\Vert u\Vert_{\alpha,\beta}=\sup_{-T<t<T}\left\vert t\right\vert
^{\frac{\alpha-\beta}{2}}\Vert u(t)\Vert_{(\rho+2,\infty)},
\]
which are weakly continuous in the sense of distributions at $t=0$.
\end{definition}

\begin{definition}
\label{d1} Let $0<T\leq\infty.$ A mild solution of the initial value problem (\ref{sc1}%
)-(\ref{sc2}) in the space $E_{\alpha,\beta}^{T}$ (respectively, in
the space $E_{\alpha}$) is a complex valued function
$u\!\in\!E_{\alpha,\beta}^{T}$ (respectively, $u\in E_{\alpha}$)
satisfying equation (\ref{sc6}) for all $0<\left\vert t\right\vert
<T$, such that $u(t)\rightharpoonup\phi\mbox{ when }t\rightarrow 0$
in the sense of distributions.
\end{definition}

Our main results are

\begin{theorem}
\label{teogeral} (Local in time solutions) Let $0<\rho<\infty$ and
$\displaystyle \frac{n\rho}{2}<\frac{\rho+2}{\rho+1}$.

\begin{enumerate}
\item If $\phi\in L^{(\frac{\rho+2}{\rho+1},\infty)}$, then there exists
$0<T<\infty$ such that the initial value problem (\ref{sc1})-(\ref{sc2}) has a
unique mild solution $u(t,x)\in E_{\alpha,\beta}^{T}$, with $T=T(\phi
)=C\left\Vert \phi\right\Vert _{(\frac{\rho+2}{\rho+1},\infty)}^{-\frac{\rho
}{\delta}},$ where $\delta=1-\frac{\alpha-\beta}{2}(\rho+1)>0$.

\item Moreover, if $\phi_{n}\in L^{(\frac{\rho+2}{\rho+1},\infty)}$ is a
sequence of functions satisfying $\phi_{n}\rightarrow\phi$ in $L^{(\frac
{\rho+2}{\rho+1},\infty)},$ then there exists $0<T_{0\text{ }}<\infty$ and
$n_{0}\in\mathbb{N}$ such that, for $n\geq n_{0}$, the solutions $u_{n}$ and
$u$ with respective initial data $\phi_{n}$ and $\phi$ lie in $E_{\alpha
,\beta}^{T_{0}}$ and $u_{n}\rightarrow u$ in $E_{\alpha,\beta}^{T_{0}}$.
Actually, the solution map $\phi\mapsto u$ is Lipschitz continuous.
\end{enumerate}
\end{theorem}

\begin{theorem}
\label{teo2}(Global in time solutions) Let $0<\rho<\infty$ and
$\displaystyle \frac{\rho+2}{\rho+1}<\frac{n\rho}{2}<\rho+2.$

\begin{enumerate}
\item If $\phi$ is a distribution such that $\sup_{-\infty<t<\infty}%
|t|^{\frac{\alpha}{2}}\Vert S(t)\phi\Vert_{(\rho+2,\infty)}<\varepsilon$, for
$\varepsilon>0$ small enough, then the initial value problem (\ref{sc1}%
)-(\ref{sc2}) has a global in time mild solution $u(t,x)\in E_{\alpha}$. This
solution is the only one satisfying $\Vert u\Vert_{\alpha}\leq2\varepsilon$.

\item Futhermore, if $(\phi_{n})$ is a sequence of distributions such that
$\Vert S(t)\phi_{n}-S(t)\phi\Vert_{E_{\alpha}}$ $\rightarrow0$ when
$n\rightarrow\infty$, and $u_{n},u$ are the solutions with respective initial
data $\phi_{n}$ and $\phi$, then $u_{n}\rightarrow u$ in $E_{\alpha}.$
\end{enumerate}
\end{theorem}

\label{RemarkTeo} We compare the theorems above with previous results.

\begin{itemize}
\item In \cite{Cazenave5}, the existence of solutions in spaces of infinite
energy was obtained for $\rho_{0}<\rho<\frac{4}{n-2}$, where $\rho_{0}$ is the
value of $\rho$ for which $\frac{\rho+2}{\rho+1}=\frac{n\rho}{2}$. In
\cite{Cazenave-Vega}, using Strichartz-type inequalities, the existence of
global solutions in the class $L^{(p,\infty)}(R^{n+1})$ $\equiv L_{t}%
^{(p,\infty)}\left(  L_{x}^{(p,\infty)}\right) $ was established, where
$p=\frac{\rho(n+2)}{2(\rho+1)}$ and $\rho_{0}<$ $\frac{4(n+1)}{n(n+2)}%
<\rho<\frac{4(n+1)}{n^{2}}.$
%On the other hand, in \cite{Planchon1}, using an
%different approach based in Besov spaces, the author has analyzed the
%existence of solutions in the Besov space $\dot{B}_{2}^{s_{\rho},\infty}%
%(R^{n})$ with positive regularity $s_{\rho}=\frac{n}{2}-\frac{2}{\rho}>0$ and
%$\rho_{0}<\frac{4}{n}<\rho<\infty.$
So, Theorem \ref{teogeral} extends the set of exponents $\rho$ where such
solutions exist by including the interval $0<\rho<\rho_{0}$.

\item In the range $\rho_{0}<\rho<\frac{4}{n-2}$, Theorem \ref{teo2} extends
the global solutions results derived in \cite{Cazenave5} to the framework of
Marcinkiewicz spaces. Our range for $\rho$ is also greater than the one in
\cite{Cazenave-Vega} (see \ref{range2}).

\item Theorem \ref{teogeral} assures the existence of local in time solutions
even for singular initial data $\phi(x)=P_{k}(x)\left\vert
x\right\vert ^{-k-\frac{n(\rho+1)}{\rho+2}}\in
L^{(\frac{\rho+2}{\rho+1},\infty)}$, where $P_{k}(x)$ is a
homogeneous polynomial of degree $k$. As far as we know, there were
no previous existence results covering this case. On the other hand,
we were not able to obtain self-similar solutions in
$E_{\alpha,\beta}$ though, since the norm
$\Vert\cdot\Vert_{\alpha,\beta}$ is not invariant by the scaling
relation $u_{\mu}(t,x)=\mu^{\frac{2}{\rho}}u(\mu^{2}t,\mu x).$
\end{itemize}

%\end{remark}
As a direct consequence of Theorem \ref{teo2}, one can show the existence of a
self-similar solution.
%taking the initial data with the
%correct homogeneity, we can show the existence of self-similar solution. This
%is the content of the following corollary.

\begin{corollary}
\label{selfsim}(self-similar solutions) In addition to the hypothesis of
Theorem \ref{teo2}, if the initial data $\phi$ is a sufficiently small
homogeneous function of degree $-\frac{2}{\rho}$, then the solution $u(t,x)$
provided by Theorem \ref{teo2} is self-similar, that is, $u(t,x)=\mu^{\frac
{2}{\rho}}u(\mu^{2}t,\mu x)$ for all $\mu>0$, almost everywhere for
$x\in\mathbb{R}^{n}$ and $t>0$.
\end{corollary}

\begin{remark}
Let $P_{k}(x)$ be a homogeneous polynomial of degree $k$. The set of functions
$\phi$ which are finite linear combinations of functions of the form
$\frac{P_{k}(x)}{\left\vert x\right\vert ^{k+\frac{2}{\rho}}}$ is an
admissible class for the existence of self-similar solutions for problem
(\ref{sc1})-(\ref{sc2}).
\end{remark}

We also analyze the large time behaviour of the solutions given by Theorem
\ref{teo2}, and study the behaviour of the solutions given in Theorem
\ref{teogeral} near to time $t=0$. These are the content of the following theorem.

\begin{theorem}
\label{TeoAssin}

\begin{enumerate}
\item {(Asymptotic stability) Suppose $0\leq h<1-\frac{\alpha}{2}(\rho+1)$,
and let $u$, $v\in E_{\alpha}$ be two global solutions of problem
(\ref{sc1})-(\ref{sc2}) obtained through Theorem \ref{teo2},
corresponding to respective initial conditions $\phi$, $\varphi$
$\in L^{(\frac{\rho+2}{\rho +1},\infty)}$. If
$\displaystyle\lim_{\left\vert
t\right\vert\rightarrow\infty}\left\vert t\right\vert
^{\frac{\alpha}{2}+h}\left\Vert S(t)(\phi-\varphi)\right\Vert
_{(\rho+2,\infty)}=0$, then
\begin{equation}
\lim_{\left\vert
t\right\vert\rightarrow\infty}\left\vert t\right\vert ^{\frac{\alpha}{2}%
+h}\left\Vert u(t)-v(t)\right\Vert _{(\rho+2,\infty)}=0. \label{as1}%
\end{equation}
}

\item (Decay rate as $t\rightarrow0$) Suppose $\delta=1-\frac{\alpha-\beta}%
{2}(\rho+1)>0$, and $h>-\delta$ . Let $u,v\in E_{\alpha,\beta}$ be two local
solutions of (\ref{sc1})-(\ref{sc2}) obtained through Theorem \ref{teogeral},
corresponding to initial conditions $\phi,\varphi$ $\in L^{(\frac{\rho+2}%
{\rho+1},\infty)}$, respectively. If $\displaystyle\lim_{t\rightarrow
0}\left\vert t\right\vert ^{\frac{\alpha-\beta}{2}-h}\left\Vert S(t)(\phi
-\varphi)\right\Vert _{(\rho+2,\infty)}=0$, then
\begin{equation}
\lim_{t\rightarrow0}\left\vert t\right\vert ^{\frac{\alpha-\beta}{2}%
-h}\left\Vert u(t)-v(t)\right\Vert _{(\rho+2,\infty)}=0.\label{as2}%
\end{equation}

\end{enumerate}
\end{theorem}

Let us comment some improvements produced by Theorem \ref{TeoAssin}.

\begin{itemize}
\item (Asymptotic stability) Theorem \ref{teo2} already gives
\[
\sup_{\left\vert t\right\vert>0}\left\vert t\right\vert
^{\frac{\alpha}{2}}\left\Vert u(t)-v(t)\right\Vert
_{(\rho+2,\infty)}<\infty.
\]
Thus, it is obvious that the estimate (\ref{as1}) holds for $h<0$. On the
other hand, the first item in Theorem \ref{TeoAssin} extends this property for
the range $0\leq h<1-\frac{\alpha}{2}(\rho+1)$. However, more regularity on
the initial perturbation $\phi-\varphi$ is required though. For instance,
assuming (in addition) that $\phi-\varphi\in L^{\frac{\rho+2}{\rho+1}}$, one
obtains
\[
\lim_{\left\vert
t\right\vert\rightarrow\infty}\left\vert t\right\vert ^{\frac{\alpha}{2}%
+h}\left\Vert S(t)(\phi-\varphi)\right\Vert _{(\rho+2,\infty)}=0,
\]
with $0\leq h<-\frac{\beta}{2}$. Observe that $-\frac{\beta}{2}=1-\frac
{\alpha}{2}(\rho+1)>0,$ when $\rho_{0}<\rho<\frac{4}{n-2}.$

\item (Decay rate when $t\rightarrow0$) By bound (\ref{sc9}), one
can see that
\[
\left\vert t\right\vert ^{\frac{\alpha-\beta}{2}-h}\left\Vert
u(t)-v(t)\right\Vert _{(\rho+2,\infty)}\leq\left\vert t\right\vert
^{\frac{\alpha-\beta}{2}-h}\left\Vert S(t)(\phi-\varphi)\right\Vert
_{(\rho+2,\infty)}+C\left\vert t\right\vert ^{\delta-h},
\]
which implies the bound (\ref{as2}) for $h<\delta$. Assuming further
regularity for $\phi-\varphi$, the second item of Theorem \ref{TeoAssin}
extends this property for the range $h>-\delta$.
\end{itemize}

\section{Proofs}

\label{proofs}

The following Lemma is important to our ends. For its proof, see
\cite{Fer-Vill}.

\begin{lemma}
\label{genlem} Let $0<\rho<\infty$ and $X$ to be a Banach space with norm
$\Vert\cdot\Vert$. Suppose $B:X\rightarrow X$ to be a map satisfying
\begin{equation}
\Vert B(x)-B(z)\Vert\leq K\Vert x-z\Vert\left(  \Vert x\Vert^{\rho}+\Vert
z\Vert^{\rho}\right),  \label{p2}%
\end{equation}
$B(0)=0$, and let $R>0$ be the unique positive root of equation
$2^{\rho+1}K(R)^{\rho }-1=0$. Given $0<\varepsilon<R$ and $y\in X$,
$y\neq0$, such that $\Vert y\Vert\leq\varepsilon$, there exists a
solution $x\in X$ for the equation $x=y+B(x)$ such that $\Vert
x\Vert\leq2\varepsilon$. The solution $x$ is unique in the
ball~$B_{2\varepsilon}:=\overline{B}(0,2\varepsilon).$ Moreover, the
solution depends continuously on $y$ in the following sense: If
$\Vert\tilde{y}\Vert\leq\varepsilon$,
$\tilde{x}=\tilde{y}+B(\tilde{x})$, and
$\Vert\tilde{x}\Vert\leq2\varepsilon$, then%
\[
\Vert x-\tilde{x}\Vert\leq\frac{1}{1-2^{\rho+1}K\varepsilon^{\rho}}\Vert
y-\tilde{y}\Vert.
\]

\end{lemma}

Now, we state and prove the necessary estimates in order to apply Lemma
\ref{genlem} in our case.

\begin{lemma}
\label{lem8} Let $0<\rho<\infty$ and $B$ be defined as
\[
B(u) = -i\lambda\int_{0}^{t} S(t-s)(|u(s)|^{\rho}u(s))ds .
\]
If $\displaystyle \frac{n\rho}{2}<\frac{\rho+2}{\rho+1}$, then there exists a
positive constant $K_{\alpha,\beta}$ such that%
\begin{equation}
\Vert B(u)-B(v)\Vert_{\alpha,\beta}\leq{K}_{\alpha,\beta}T^{1-\frac
{(\alpha-\beta)(\rho+1)}{2}}\Vert u-v\Vert_{\alpha,\beta}\left(  \Vert
u\Vert_{\alpha,\beta}^{\rho}+\Vert v\Vert_{\alpha,\beta}^{\rho}\right)  ,
\label{sc9}%
\end{equation}
for all $u$, $v \in \!E_{\alpha,\beta}^{T}$. On the other hand, if
$\displaystyle \frac{\rho+2}{\rho+1}<\frac{n\rho}{2}<\rho+2$, then
there exists a positive constant $K_{\alpha}$ such that
\begin{equation}
\Vert B(u)-B(v)\Vert_{\alpha}\leq{K_{\alpha}} \Vert u-v\Vert_{\alpha}\left(
\Vert u\Vert_{\alpha}^{\rho}+\Vert v \Vert_{\alpha}^{\rho}\right)  ,
\label{sc8}%
\end{equation}
for all $u$, $v\in E_{\alpha}$.
\end{lemma}

\begin{proof} Without loss of generality, we assume $t>0$.
First note that if $\frac{n\rho}{2}<\frac{\rho+2}{\rho+1}<\rho+2$, then
$\frac{\alpha-\beta}{2}(\rho+1)<1$ and $\frac{n}{2}(\frac{2(\rho+1)}{\rho
+2}-1)<1$. Therefore,
\begin{align*}
\Vert B(u)-B(v)\Vert_{({\rho+2,\infty)}}  &  \leq\int_{0}^{t}\Vert
S(t-s)(\left\vert u\right\vert ^{\rho}u-\left\vert v\right\vert ^{\rho}%
v)\Vert_{({\rho+2,\infty)}}ds\\
&  \leq C\int_{0}^{t}(t-s)^{-\frac{n}{2}(\frac{2(\rho+1)}{\rho+2}-1)}%
\Vert(\left\vert u-v\right\vert )(\left\vert u\right\vert ^{\rho}+\left\vert
v\right\vert ^{\rho})\Vert_{(\frac{{\rho+2}}{\rho+1}{,\infty)}}ds\\
&  \leq C\int_{0}^{t}(t-s)^{-\frac{n}{2}(\frac{2(\rho+1)}{\rho+2}-1)}\Vert
u-v\Vert_{(\rho+2{,\infty)}}\left(  \Vert u\Vert_{(\rho+2{,\infty)}}^{\rho
}+\Vert v\Vert_{(\rho+2{,\infty)}}^{\rho}\right)  ds\\
&  \hspace{-4cm}\leq C
\left(  \sup_{0<t<T}t^{\frac{\alpha
-\beta}{2}}\Vert u-v\Vert_{(\rho+2{,\infty)}}\sup_{0<t<T}\left(
t^{\frac{(\alpha-\beta)\rho}{2}}\Vert u\Vert_{(\rho+2{,\infty)}}^{\rho
}+t^{\frac{(\alpha-\beta)\rho}{2}}\Vert v\Vert_{(\rho+2{,\infty)}}^{\rho
}\right)  \right)
\int_{0}^{t}(t-s)^{-\frac{\alpha-\beta}{2}}
s^{-\frac{\alpha-\beta}{2}(\rho+1)}ds  \\
&  \hspace{-4cm}=K_{\alpha,\beta}t^{-\frac{\alpha-\beta}{2}}t^{1-\frac{\alpha-\beta}{2}%
(\rho+1)}\Vert u-v\Vert_{\alpha,\beta}\left(  \Vert
u\Vert_{\alpha,\beta }^{\rho}+\Vert
v\Vert_{\alpha,\beta}^{\rho}\right)  ,
\end{align*}
which proves (\ref{sc9}). On the other hand, if
$\frac{\rho+2}{\rho+1}<$ $\frac{n\rho}{2}<\rho+2,$ then
$\frac{\alpha}{2}(\rho+1)<1$ and
$\frac{n}{2}(\frac{2(\rho+1)}{\rho+2}-1)<1$. In this case,
\begin{align*}
\Vert B(u)-B(v)\Vert_{({\rho+2,\infty)}}  &  \leq C\int_{0}^{t}(t-s)^{-\frac
{n}{2}(\frac{2(\rho+1)}{\rho+2}-1)}\Vert u-v\Vert_{(\rho+2{,\infty)}}\left(
\Vert u\Vert_{(\rho+2{,\infty)}}^{\rho}+\Vert v\Vert_{(\rho+2{,\infty)}}%
^{\rho}\right)  ds\\
&  \hspace{-3.2cm}\leq C \left(  \sup_{t>0}t^{\frac{\alpha
}{2}}\Vert u-v\Vert_{(\rho+2{,\infty)}}\sup_{t>0}\left(  t^{\frac{\alpha\rho
}{2}}\Vert u\Vert_{(\rho+2{,\infty)}}^{\rho}+t^{\frac{\alpha\rho}{2}}\Vert
v\Vert_{(\rho+2{,\infty)}}^{\rho}\right)  \right)\int_{0}^{t}(t-s)^{-\frac{n}{2}(\frac{2(\rho+1)}
{\rho+2}-1)}s^{-\frac{\alpha}{2}(\rho+1)}ds \\
&  \hspace{-3.2cm}=K_{\alpha}t^{-\frac{\alpha}{2}}\Vert
u-v\Vert_{\alpha}\left(  \Vert u\Vert_{\alpha}^{\rho}+\Vert
v\Vert_{\alpha}^{\rho}\right)  ,
\end{align*}
which proves (\ref{sc8}). \end{proof}

\subsection{Proof of Theorem \ref{teogeral}}

%We apply Lemma \ref{genlem} to the integral equation (\ref{sc6}) with
%$X=E_{\alpha,\beta}^{T}$ and $y=S(t)\phi.$
Let $y=S(t)\phi$. Due to Lemma \ref{estim1}, one has
\[
\left\Vert y\right\Vert _{\alpha,\beta}=\sup_{-T<t<T}|t|^{\frac{^{\alpha
-\beta}}{2}}\left\Vert S(t)\phi\right\Vert _{(\rho+2,\infty)}\leq C\left\Vert
\phi\right\Vert _{(\frac{\rho+2}{\rho+1},\infty)}<\infty.
\]
Using Lemma \ref{lem8}, one gets
\begin{equation}
\Vert B(u)-B(v)\Vert_{\alpha,\beta}\leq
K_{\alpha,\beta}T^{\delta}\Vert u-v\Vert _{\alpha,\beta}\left(
\Vert u\Vert_{\alpha,\beta}^{\rho}+\Vert v\Vert
_{\alpha,\beta}^{\rho}\right)  , \label{sc10}%
\end{equation}
where $\delta=1-\frac{\alpha-\beta}{2}(\rho+1)>0.$ Now, choose
$0<T<\infty$ sufficiently small, and $\varepsilon>0$ such that
$\left\Vert y\right\Vert _{\alpha,\beta}\leq C\left\Vert
\phi\right\Vert _{(\frac{\rho+2}{\rho
+1},\infty)}=\varepsilon<R:=\left(
\frac{1}{2^{(\rho+1)}K_{\alpha,\beta}T^{\delta}}\right)
^{\frac{1}{\rho}}$. Using Lemma \ref{genlem} with
$X=E_{\alpha,\beta}^{T}$,
one assures the existence of a local mild solution $u\in E_{\alpha,\beta}^{T}%
$. Moreover, this solution is unique in the ball ${B}_{2\varepsilon
}:=\overline{B}(0,2\varepsilon)\subset E_{\alpha,\beta}^{T}$.
Furthermore, through standard arguments on can prove that
$u(t)\rightarrow\phi$ in the sense of distributions when
$t\rightarrow0$. So, solutions of the integral equation are indeed
mild solutions in the sense of Definition \ref{d1}.

Finally, let $u_{n}$ and $u$ be the solutions with respective initial data
$\phi_{n}$ and $\phi$. By Lemma \ref{genlem}, one has
\[
\Vert
u_{n}-u\Vert_{E_{\alpha,\beta}}\leq\frac{1}{1-2^{\rho+1}K_{\alpha,\beta}T^{\delta}\varepsilon
^{\rho}}\Vert S(t)\phi_{n}-S(t)\phi\Vert_{E_{\alpha,\beta}}\leq\frac
{C}{1-2^{\rho+1}K_{\alpha,\beta}T^{\delta}\varepsilon^{\rho}}\Vert\phi_{n}-\phi\Vert_{(\frac{\rho
+2}{\rho+1},\infty)}.
\]
This finishes the proof.

\subsection{Proof of Theorem \ref{teo2}}

Apply Lemma \ref{genlem} to the integral equation (\ref{sc6}) with
$X=E_{\alpha}$ and $y=S(t)\phi$. In this case, the bound (\ref{sc8}) gives
\[
\Vert B(u)-B(v)\Vert_{\alpha}\leq K_{\alpha}\Vert
u-v\Vert_{\alpha}\left( \Vert u\Vert_{\alpha}^{\rho}+\Vert
v\Vert_{\alpha}^{\rho}\right)  .
\]
Now, one considers $\varepsilon> 0$ small enough so that $\Vert
S(t)\phi \Vert_{\alpha}=\sup_{\left\vert
t\right\vert>0}|t|^{\frac{\alpha}{2}}\Vert S(t)\phi\Vert
_{(\rho+2,\infty)} < \varepsilon$ allows one to apply Lemma
\ref{genlem} repeatedly, in order to obtain the existence of a
global mild solution $u\in E_{\alpha}$. This solution is unique in
the ball ${B}_{2\varepsilon }:=\overline{B}(0,2\varepsilon)\subset
E_{\alpha}$.

The continuity of the solutions with respect to the initial conditions, as
well as the continuity of the solutions in the sense of distributions, follow
as in the proof of Theorem \ref{teogeral}.

\subsection{Proof of Corollary \ref{selfsim}}

%If $u(t,x)$ is a smooth solution of the Schr\"{o}dinger equation
%(\ref{sc1})-(\ref{sc2}), then it is straightforward to check that for all
%$\mu>0$, the function $u_{\mu}(t,x)=\mu^{\frac{2}{\rho}}u(\mu^{2}t,\mu x)$ is
%also solves (\ref{sc1})-(\ref{sc2}). Solutions of (\ref{sc1})-(\ref{sc2})
%satisfying
%\begin{equation}
%u(t,x)=u_{\mu}(t,x){,} \label{autosimi1}%
%\end{equation}
%for any $\mu>0,$ $t\in\mathbb{R},$ and $x\in\mathbb{R}^{n}$, are the so called
%self-similar solutions. In order to prove the existence of such solutions,
Let $t>0$. If the initial data $\phi(x)$ is a homogeneous function
of degree $-\frac {2}{\rho}$, then $S(t)\phi$ satisfies the
self-similar property $u(t,x) = \mu^{\frac{2}{\rho}}u(\mu^{2}t,\mu
x)$. Thus,
\[
t^{\alpha/2}\Vert S(t)\phi\Vert_{(\rho+2,\infty)}=t^{\frac{\alpha}{2}}%
t^{\frac{n}{2(\rho+2)}-\frac{1}{\rho}}\Vert S(1)\phi\Vert_{(\rho+2,\infty
)}=\Vert S(1)\phi\Vert_{(\rho+2,\infty)}.
\]
Moreover, $\Vert S(1)\phi\Vert_{L^{\rho+2}}$ is finite(see \cite{Cazenave5}).
Since the inclusion $L^{\rho+2}\hookrightarrow L^{(\rho+2,\infty)}$
continuously, one has
\[
\Vert S(1)\phi\Vert_{(\rho+2,\infty)}\leq\Vert S(1)\phi\Vert_{L^{\rho+2}%
}<\infty.
\]
Therefore, if $\Vert S(1)\phi\Vert_{(\rho+2,\infty)}$ is small enough , it is
straightforward to show that the solution $u(t,x)$ obtained in Theorem
\ref{teogeral} is self-similar.

\subsection{Proof of Theorem \ref{TeoAssin}}

Without loss of generality, assume $t>0$. Subtracting the integral
equations satisfied by $u$ and $v$, one gets
\begin{align*}
t^{\frac{\alpha}{2}+h}\left\Vert u(t)-v(t)\right\Vert _{(\rho+2,\infty)}  &
\leq t^{\frac{\alpha}{2}+h}\left\Vert S(t)(\phi-\varphi)\right\Vert
_{(\rho+2,\infty)}\\
& \mbox{}+t^{\frac{\alpha}{2}+h}\left\Vert \int_{0}^{t}S(t-s)(u\left\vert
u\right\vert ^{\rho}-v\left\vert v\right\vert ^{\rho})ds\right\Vert
_{(\rho+2,\infty)}.
\end{align*}
Since $\Vert u\Vert_{\alpha},\Vert v\Vert_{\alpha}\leq2\varepsilon$, one uses
the change of variable $s\longmapsto ts$, and bound
\begin{align*}
&  t^{\frac{\alpha}{2}+h}\Vert\int_{0}^{t}S(t-s)(u|u|^{\rho}-v|v|^{\rho
})ds\Vert_{(\rho+2,\infty)}\\
&  \leq Ct^{\frac{\alpha}{2}+h}\int_{0}^{t}(t-s)^{-\frac{n}{2}(\frac
{2(\rho+1)}{(\rho+2)}-1)}s^{-\frac{\alpha(\rho+1)}{2}-h}(s^{\frac{\alpha\rho
}{2}}\Vert u(s)\Vert_{(\rho+2,\infty)}^{\rho}+s^{\frac{\alpha\rho}{2}}\Vert
v(s)\Vert_{(\rho+2,\infty)}^{\rho})s^{\frac{\alpha}{2}+h}\Vert u(s)-v(s)\Vert
_{(\rho+2,\infty)}ds\\
&  \leq C2^{\rho+1}\varepsilon^{\rho}\int_{0}^{1}(1-s)^{-\frac{n}{2}%
(\frac{2(\rho+1)}{(\rho+2)}-1)}s^{-\frac{\alpha(\rho+1)}{2}-h}(ts)^{\frac
{\alpha}{2}+h}\Vert u(ts)-v(ts)\Vert_{(\rho+2,\infty)}ds\text{.}%
\end{align*}
Therefore,
\begin{align}
t^{\frac{\alpha}{2}+h}\left\Vert u(t)-v(t)\right\Vert _{(\rho+2,\infty)}  &
\leq t^{\frac{\alpha}{2}+h}\left\Vert S(t)(\phi-\varphi)\right\Vert
_{(\rho+2,\infty)}\label{sta3}\\
& \mbox{}+C2^{\rho+1}\varepsilon^{\rho}\int_{0}^{1}(1-s)^{-\frac{n}{2}%
(\frac{2(\rho+1)}{(\rho+2)}-1)}s^{-\frac{\alpha}{2}(\rho+1)-h}(ts)^{\frac
{\alpha}{2}+h}\Vert u(ts)-v(ts)\Vert_{(\rho+2,\infty)}ds,\nonumber
\end{align}
for all $t>0$. Now, define
\[
A:=\limsup_{t\rightarrow\infty}t^{\frac{\alpha}{2}+h}\Vert u(t)-v(t)\Vert
_{(\rho+2,\infty)}.
\]
Using the assumption on the initial perturbation $\phi-\varphi$, it
is not difficult to show that $A<\infty.$ Now, note that
\[%
\begin{array}
[c]{lr}%
\displaystyle\limsup_{t\rightarrow\infty}\int_{0}^{1}(1-s)^{-\frac{n}{2}%
(\frac{2(\rho+1)}{(\rho+2)}-1)}s^{\frac{-\alpha(\rho+1)}{2}-h}(ts)^{\frac
{\alpha}{2}+h}\Vert u(ts)-v(ts)\Vert_{(\rho+2,\infty)}ds & \vspace{0.4cm}\\
& \hspace{-5cm}\leq\displaystyle A\int_{0}^{1}(1-s)^{-\frac{n}{2}(\frac
{2(\rho+1)}{(\rho+2)}-1)}s^{-\frac{\alpha}{2}(\rho+1)-h}ds\text{.}%
\end{array}
\]
So, taking $\displaystyle\limsup_{t\rightarrow\infty}$ in \eqref{sta3}, one
obtains
\[
A\leq\left(  C2^{\rho+1}\varepsilon^{\rho}\int_{0}^{1}(1-s)^{-\frac{n}%
{2}(\frac{2(\rho+1)}{(\rho+2)}-1)}s^{-\frac{\alpha(\rho+1)}{2}-h}ds\right)  A.
\]
Now, let $\Gamma:=C2^{\rho+1}\int_{0}^{1}(1-s)^{-\frac{n}{2}(\frac{2(\rho
+1)}{(\rho+2)}-1)}s^{-\frac{\alpha(\rho+1)}{2}-h}ds$. Choosing $\varepsilon>0$
sufficiently small such that $\varepsilon^{\rho}\Gamma<1$, one concludes that
$A=0$. This proves part 1 of the theorem.

In order to prove part 2, let
$\delta=1-\frac{\alpha-\beta}{2}(\rho+1)$ and $0<t<T$  as in Theorem
\ref{teogeral}. One can write $\delta=\frac
{\alpha-\beta}{2}-h-\frac{\alpha-\beta}{2}-\frac{\alpha-\beta}{2}%
(\rho+1)+h+1.$ Again, one subtracts the equations for $u$ and $v$, and bound
\begin{align*}
&  t^{\frac{\alpha-\beta}{2}-h}\left\Vert \int_{0}^{t}S(t-s)(u\left\vert
u\right\vert ^{\rho}-v\left\vert v\right\vert ^{\rho})ds\right\Vert
_{(\rho+2,\infty)}\\
&  \leq Ct^{\frac{\alpha-\beta}{2}-h}\int_{0}^{t}(t-s)^{-\frac{\alpha-\beta
}{2}}s^{-\frac{\alpha-\beta}{2}(\rho+1)+h}(s^{\frac{\alpha-\beta}{2}\rho
}(\Vert u(s)\Vert_{(\rho+2,\infty)}^{\rho}+\Vert v(s)\Vert_{(\rho+2,\infty
)}^{\rho}))\text{ }s^{\frac{\alpha-\beta}{2}-h}\Vert u(s)-v(s)\Vert
_{(\rho+2,\infty)}ds\\
&  \leq C2^{\rho+1}\varepsilon^{\rho}t^{\frac{\alpha-\beta}{2}-h-\frac
{\alpha-\beta}{2}-\frac{\alpha-\beta}{2}(\rho+1)+h+1}\int_{0}^{1}%
(1-s)^{-\frac{\alpha-\beta}{2}}s^{-\frac{\alpha-\beta}{2}(\rho+1)+h}%
(ts)^{\frac{\alpha-\beta}{2}-h}\Vert u(ts)-v(ts)\Vert_{(\rho+2,\infty)}ds.
\end{align*}
Hence,
\begin{align}
t^{\frac{\alpha-\beta}{2}-h}\left\Vert u(t)-v(t)\right\Vert _{(\rho
+2,\infty)}  &  \leq t^{\frac{\alpha-\beta}{2}-h}\left\Vert S(t)(\phi
-\varphi)\right\Vert _{(\rho+2,\infty)}\nonumber\\
&  \mbox{}+C2^{\rho+1}\varepsilon^{\rho}t^{\delta}\int_{0}^{1}(1-s)^{-\frac
{\alpha-\beta}{2}}s^{-\frac{\alpha-\beta}{2}(\rho+1)+h}(ts)^{\frac
{\alpha-\beta}{2}-h}\Vert u(ts)-v(ts)\Vert_{(\rho+2,\infty)}ds\nonumber
\end{align}
Writing
$A:=\lim\sup_{t\rightarrow0}t^{\frac{\alpha-\beta}{2}-h}\Vert
u(t)-v(t)\Vert_{(\rho+2,\infty)}<\infty$, one takes $\displaystyle
\lim \sup_{t\rightarrow0}$ in the last inequality to obtain
\[
0\leq A\leq C2^{\rho+1}\varepsilon^{\rho}A\int_{0}^{1}(1-s)^{-\frac{\alpha-\beta}%
{2}}s^{-\frac{\alpha-\beta}{2}(\rho+1)+h}ds\text{
}\lim_{t\rightarrow 0}t^{\delta}=0.
\]This concludes the proof.

\end{document}